\documentclass[10pt]{article}
\usepackage{}
\usepackage{amssymb}
\usepackage{amsfonts}
\usepackage{bbm}
\usepackage{mathrsfs}
\usepackage{mathrsfs}
\usepackage{longtable}
\usepackage{supertabular}
\usepackage{amsfonts,mathrsfs,latexsym,amsmath,amssymb,amsthm}

\newtheorem{proposition}{Proposition}
\newtheorem{lemma}{Lemma}

\begin{document}
\title{\bf{On generalized quasi-cyclic codes over $\mathbb{Z}_4$}}
\author{{\bf  Jian Gao$^1$, Xiangrui Meng$^1$, Fang-Wei Fu$^2$}\\
 {\footnotesize \emph{ 1. Shandong University of Technology}}\\
  {\footnotesize  \emph{Zibo, 255091, P. R. China}}\\
 {\footnotesize \emph{ 2. Chern Institute of Mathematics and LPMC, Nankai University}}\\
  {\footnotesize  \emph{Tianjin, 300071, P. R. China}}\\
 }
\date{}

\maketitle \noindent {\small {\bf Abstract} Based on good algebraic structures and practicabilities, generalized quasi-cyclic (GQC) codes play important role in coding theory. In this paper, we study some results on GQC codes over $\mathbb{Z}_4$ including the normalized generating set, the minimum generating set and the normalized generating set of their dual codes. As an application, new $\mathbb{Z}_4$-linear codes and good nonlinear binary codes are constructed from GQC codes over $\mathbb{Z}_4$.}
 \vskip 1mm

\noindent
 {\small {\bf Keywords}  Generalized quasi-cyclic codes; Normalized generating sets; Minimum generating sets; New $\mathbb{Z}_4$-linear codes}

\vskip 3mm \noindent {\bf Mathematics Subject Classification (2000) } 11T71 $\cdot$ 94B05 $\cdot$ 94B15

\vskip 3mm \baselineskip 0.2in

\section{Introduction}
\par
The study of codes over rings began in 1970s. In 1994, Hammons et al. studied $\mathbb{Z}_4$-linear codes extensively, and they showed that some binary nonlinear codes with good parameters can be viewed as Gray images of some cyclic codes over $\mathbb{Z}_4$ \cite{Hammons}. From then on, studying codes over rings, especially codes over the quaternary ring $\mathbb{Z}_4$, has become a hot research topic. Many classes of $\mathbb{Z}_4$-codes such as linear cyclic codes, LCD codes, projective codes, trace codes and so on were studied by coding scientists \cite{Dougherty,Shi1,Shi2,Shi3,Shi4}.
 \par
Quasi-cyclic (QC) codes have been studied extensively in the last years. One motivation of researching QC codes is their good algebraic structure for encoding and decoding. Another motivation is that QC codes contain a lot of good linear codes over finite fields and finite rings. Further, QC codes can produce good lattices, quantum codes and convolutional codes. In \cite{Ling}, Ling and Sol\'{e} studied QC codes over finite chain rings by the Chinese remainder theorem. They gave a new quaternary construction of Lee lattices. By the Gray map, new binary linear codes are obtained from QC codes with odd length components over $\mathbb{Z}_4$ \cite{Aydin}. In \cite{Siap1}, QC codes with even length components over $\mathbb{Z}_4$ were studied by Siap et al.. They determined the minimum generating set of QC codes, and obtained new binary nonlinear codes from the usual Gray map from $\mathbb{Z}_4$ to $\mathbb{Z}_2$.
\par
As the generalization of QC codes, generalized quasi-cyclic (GQC) codes were studied widely in recent years. Basic algebraic structures of GQC codes were firstly introduced in \cite{Siap}. Afterwards, Esmaeili and Yari studied GQC codes extensively by the Chinese remainder theorem \cite{Esmaeili}. In \cite{Guneri2}, G\"{u}neri et al. decomposed GQC codes into a direct sum of linear codes and gave a lower bound on the minimum Hamming distance, the trace representation of GQC codes. They also showed that GQC codes are asymptotically good. Recently, the concept of GQC codes has been generalized to finite rings. In \cite{Cao}, Cao studied the enumeration of GQC codes over Galois rings. In \cite{Gao1}, we studied GQC codes over the ring $\mathbb{F}_q+u\mathbb{F}_q$ including minimum generating sets, minimum distance and codes construction. For GQC codes with only one given formal generator over $\mathbb{Z}_4$, we studied their the minimum generating set, the minimum Lee distance and constructed good binary nonlinear codes \cite{Wu}. Not that, in \cite{Wu}, we did not give the explicit generators of GQC codes.
\par
Another point of view for researching GQC codes is that how to determine their explicit generators, which similar to the generator polynomial of cyclic codes are important to determining the generator matrix, dual codes and the minimum distance of GQC codes. In \cite{Borges2}, Borges et al. firstly determined the explicit generator for $\mathbb{Z}_2$-double cyclic codes, which are GQC codes with index 2 over $\mathbb{Z}_2$. Following \cite{Borges2}, we determined the explicit generator for $\mathbb{Z}_4$-double cyclic codes, and gave the relationship between $\mathbb{Z}_4$-double cyclic codes and their dual codes  \cite{Gao}. As an application, using the Gray map, some good binary nonlinear codes were obtained. Further, we proved that $\mathbb{Z}_4$-double cyclic codes are asymptotically good \cite{Gao2}. Recently, Bae et al. studied the explicit generator for GQC codes of arbitrary index over $\mathbb{Z}_2$ in \cite{Bae}.
\par
In this paper, motivated by the above papers, we will determine explicit generators of GQC codes and their dual codes with the arbitrary index over $\mathbb{Z}_4$. Our main contributions in this paper are listed below.
\par
(i) This paper determined the normalized generating set and the minimum generating set of GQC codes with index $l$ over $\mathbb{Z}_4$, where $l$ is a positive integer.
\par
(ii) This paper determined the dual codes of GQC codes over $\mathbb{Z}_4$.
\par
The rest of this paper is organized as follows. In Section 2, we review some results on polynomials and codes over $\mathbb{Z}_4$. In Section 3, we introduce some definitions and give the normalized generating set of GQC codes over $\mathbb{Z}_4$. In Section 4, we determined the minimal generating set of GQC codes over $\mathbb{Z}_4$. In Section 5, we give the relationship between GQC codes and their dual codes over $\mathbb{Z}_4$. 
\section{Preliminary}
Let $\mathbb{Z}_4=\{0,1,2,3\}$ be a quaternary ring. Any element $c\in \mathbb{Z}_4$ can be written as $c=a+2b$, where $a,b\in \mathbb{Z}_2$. Define the map $-$ from $\mathbb{Z}_4$ to $\mathbb{Z}_2$ as $\overline{a+2b}=a$. Extend it to the polynomial ring
\begin{align*}
-: \mathbb{Z}_4[x] &\rightarrow \mathbb{Z}_2[x]\\
a_0+a_1x+\cdots+a_nx^n&\mapsto \overline{a}_0+\overline{a}_1x+\cdots+\overline{a}_nx^n.
\end{align*}
Clearly, it is a surjective ring homomorphism map.
\par
 Let $f(x)$ and $g(x)$ be two polynomials over $\mathbb{Z}_4$. If there exist polynomials $u(x),v(x)\in \mathbb{Z}_4[x]$ such that
$u(x)f(x)+v(x)g(x)=1$, then $f(x)$ and $g(x)$ are coprime over $\mathbb{Z}_4$.
Further, $f(x)$ and $g(x)$ are coprime over $\mathbb{Z}_4$ if and only if $\overline{f}(x)$ and $\overline{g}(x)$ are coprime over $\mathbb{Z}_2$  (see Lemma 5.1 in \cite{Wan}).
Let $f(x)$ be a monic polynomial over $\mathbb{Z}_4$.
If $\overline{f}(x)$ is irreducible over $\mathbb{Z}_2$, then $f(x)$ is called basic irreducible over $\mathbb{Z}_4$.
\begin{lemma}\cite [Lemma 5.3]{Wan} (Hensel's Lemma)
Let $f(x)$ be a monic polynomial of degree $\geq 1$ over $\mathbb{Z}_4$. The factorization of $\overline{f}(x)$ over $\mathbb{Z}_2$ is
$$\overline{f}(x)=t_1(x)t_2(x)\cdots t_r(x),$$
where $t_1(x),t_2(x),\ldots,t_r(x)$ are pairwise coprime polynomials over $\mathbb{Z}_2$.
Then there are pairwise coprime monic polynomials $f_i(x)\in \mathbb{Z}_4[x]$ for $i=1,2,\ldots,r$ such that $f(x)=f_1(x)f_2(x)\cdots f_r(x)$, $\overline{f}_i(x)=t_i(x)$ and ${\rm deg}(f_i(x))={\rm deg}(t_i(x))$ for $i=1,2,\ldots,r$.
\end{lemma}
\par
Let $I$ be an ideal of $\mathbb{Z}_4[x]$. If $I\neq \mathbb{Z}_4[x]$ and $a(x)b(x)\in I$ implies $a(x)\in I$ or $b(x)^n\in I$ for some positive integer $n$, then the ideal $I$ is said to be primary. Let $a(x)$ be a nonzero polynomial over $\mathbb{Z}_4$. If the ideal $\langle a(x)\rangle$ of $\mathbb{Z}_4[x]$ is primary, then $a(x)$ is called a primary polynomial over $\mathbb{Z}_4$. Any basic irreducible polynomial in $\mathbb{Z}_4[x]$ is primary (see Corollary 5.9 in \cite{Wan}).
\begin{lemma}\cite [Theorem 5.10]{Wan}(Unique Factorization Theorem)
Let $f(x)$ be a monic polynomial of degree $\geq 1$ over $\mathbb{Z}_4$. Then
$f(x)=f_1(x)f_2(x)\cdots f_r(x)$, where $f_1(x), f_2(x),\ldots,f_r(x)$ are pairwise coprime monic primary polynomials. Further, $f_1(x), f_2(x), \ldots,f_r(x)$ are uniquely determined up to a rearrangement.
\end{lemma}
\par
By the Hensel's Lemma and unique factorization of monic polynomials over $\mathbb{Z}_4$, we can get the following lemma directly.
\begin{lemma}
Let $f(x)$ be a monic polynomial of degree $\geq 1$ over $\mathbb{Z}_4$ and $\overline{f}(x)$ has no multiple root. Then
$$f(x)=f_1(x)f_2(x)\cdots f_r(x),$$
where $f_1(x),f_2(x),\ldots,f_r(x)$ are pairwise coprime monic basic irreducible polynomials over $\mathbb{Z}_4$. Moreover, $f_1(x),f_2(x),\ldots,f_r(x)$ are uniquely determined up to a rearrangement.
\end{lemma}
\par
Let $f(x)$ and $g(x)$ be two monic polynomials over $\mathbb{Z}_4$, where $\overline{f}(x)$ and $\overline{g}(x)$ have no multiple roots.
By Lemma 2.3, the factorizations of $f(x)$ and $g(x)$ over $\mathbb{Z}_4$ are
$$f(x)=h_1(x)^{e_1}h_2(x)^{e_2}\cdots h_r(x)^{e_r},~
g(x)=h_1(x)^{t_1}h_2(x)^{t_2}\cdots h_r(x)^{t_r},$$
where $e_i, t_i\in \{0,1\}$ for $i=1,2,\ldots,r$, $h_1(x),h_2(x),\ldots,h_r(x)$ are all pairwise coprime monic basic irreducible divisors of $f(x)$ and $g(x)$.
Define $$d(x)=h_1(x)^{{\rm min}\{e_1,t_1\}}h_2(x)^{{\rm min}\{e_2,t_2\}}\cdots h_r(x)^{{\rm min}\{e_r,t_r\}}.$$
Clearly, $d(x)$ is a monic common divisor of $f(x)$ and $g(x)$ with the highest degree, and it is called the greatest common divisor of $f(x)$ and $g(x)$ denoted by ${\rm gcd_4}(f(x),g(x))$.
If $\frac{f(x)}{{\rm gcd}_4(f(x),g(x))}$ and $\frac{g(x)}{{\rm gcd}_4(f(x),g(x))}$ are coprime over $\mathbb{Z}_4$,
then there exist $u(x),~v(x)\in \mathbb{Z}_4[x]$ such that $u(x)\frac{f(x)}{{\rm gcd}_4(f(x),g(x))}+v(x)\frac{g(x)}{{\rm gcd}_4(f(x),g(x))}=1$,
which implies that $\langle{\rm gcd}_4(f(x), g(x))\rangle\subseteq \langle f(x), g(x)\rangle$. Moreover, we have $\langle f(x), g(x)\rangle\subseteq\langle{\rm gcd}_4(f(x), g(x))\rangle$. Hence, $\langle f(x), g(x)\rangle=\langle{\rm gcd}_4(f(x), g(x))\rangle$.
\par
In the following, as the preparation, we will introduce some basic results on quaternary codes. For more details, one can refer the reference \cite{Wan}.
\par
Let $\mathbb{Z}_4^{n}$ be the $n$-tuples over $\mathbb{Z}_4$. If $\mathcal{C}$ is a nonempty subset of $\mathbb{Z}_4^{n}$, then $\mathcal{C}$ is called a quaternary code of length $n$.
If the quaternary code $\mathcal{C}\subseteq \mathbb{Z}_4^{n}$ is a $\mathbb{Z}_4$-submodule of $\mathbb{Z}_4^n$, then $\mathcal{C}$ is called a $\mathbb{Z}_4$-linear code of length $n$.
\par
The Lee weights of elements in $\mathbb{Z}_4$ are defined respectively by $$w_{L}(0)=0,~w_{L}(1)=w_{L}(3)=1,~w_{L}(2)=2.$$
For any $\mathbf{c}=(c_0,c_1,\ldots,c_{n-1}),~\mathbf{d}=(d_0,d_1,\ldots,d_{n-1})\in \mathbb{Z}_4^{n}$, the Lee weight of $\mathbf{c}$ is $w_{L}(\mathbf{c})=\sum_{i=0}^{n-1}w_{L}(c_i)$, the Lee distance of $\mathbf{c}$ and $\mathbf{d}$ is $d_{L}(\mathbf{c},\mathbf{d})=w_{L}(\mathbf{c}-\mathbf{d})$.
 If $\mathcal{C}$ is $\mathbb{Z}_4$-linear, then its minimum Lee distance is the minimum Lee weight of nonzero codewords of $\mathcal{C}$ actually.
 \par
 The Gray map $\Phi$ from $\mathbb{Z}_4$ to $\mathbb{Z}_2^2$ is defined as $\Phi(a+2b)=(b,a+b)$, where $a,b\in \mathbb{Z}_2$. Clearly, $w_L(a+2b)=w_H((b,a+b))$, where $w_H$ denotes the Hamming weight.
 For any $\mathbf{c}=(c_0,c_1,\ldots,c_{n-1})\in \mathbb{Z}_4^{n}$, define
 $$\Phi(\mathbf{c})=(b_0,b_1,\ldots,b_{n-1},a_0+b_0,a_1+b_1,\ldots,a_{n-1}+b_{n-1}).$$
If $\mathcal{C}$ is a $\mathbb{Z}_4$-code of length $n$, then $\Phi(\mathcal{C})$ is a binary code of length $2n$. Moreover, the minimum Lee weight and the minimum Lee distance of $\mathcal{C}$ are equal to the minimum Hamming weight and the minimum Hamming distance of $\Phi(\mathcal{C})$, respectively.
\section{Normalized generating set of GQC codes}
\par
For $i=1,2,\ldots,l$, let $m_i$ be odd positive integer and $n=m_1+m_2+\cdots+m_l$.
Let $\mathcal{C}$ be a $\mathbb{Z}_4$-linear code of length $n$. If for any codeword
$(c_{1,0},c_{1,1},\ldots,c_{1,m_1-1}|\cdots|c_{l,0},c_{l,1},\ldots,c_{l,m_l-1})\in \mathcal{C}$
we have
$$(c_{1,m_1-1},c_{1,0},\ldots,c_{1,m_1-2}|\cdots|c_{l,m_l-1},c_{l,0},\ldots,c_{l,m_l-2})\in \mathcal{C},$$
then $\mathcal{C}$ is called a generalized quasi-cyclic (GQC) code of block lengths $(m_1,m_2,\ldots,m_l)$ with index $l$ over $\mathbb{Z}_4$. If $l=1$, then $\mathcal{C}$ is called a $\mathbb{Z}_4$-cyclic code of length $m_1$. If $m=m_1=\cdots=m_l$, then $\mathcal{C}$ is called a quasi-cyclic (QC) code of length $ml$ with index $l$ over $\mathbb{Z}_4$.
\par
Let $R_i=\mathbb{Z}_4[x]/\langle x^{m_i}-1\rangle$ and $R=R_1\times \cdots\times R_l$. Define a $\mathbb{Z}_4$-module isomorphism map $\phi$ from $\mathbb{Z}_4^{m_1}\times \cdots\times \mathbb{Z}_4^{m_l}$ to $R$ such that
\begin{align*}
\phi((c_{1,0},c_{1,1},\ldots,c_{1,m_1-1}|\cdots|c_{l,0},c_{l,1},\ldots,c_{l,m_l-1}))
=(c_1(x)|\cdots|c_{l}(x)),
\end{align*}
where $c_i(x)=\sum_{j=0}^{m_i-1}c_{i,j}x^j$.
Define a multiplication
\begin{align*}
\ast: \mathbb{Z}_4[x]\times R&\rightarrow R\\
\alpha(x)*(c_1(x)|\cdots|c_l(x))&=(\alpha(x)c_1(x)|\cdots|\alpha(x)c_l(x)),
\end{align*}
where $\alpha(x)\in \mathbb{Z}_4[x]$ and $(c_1(x)|\cdots|c_l(x))\in R$. Clearly, $R$ is a $\mathbb{Z}_4[x]$-module under this multiplication and the usual addition of vectors,
which implies that $\mathcal{C}$ is a GQC code if and only if $\phi(\mathcal{C})$ is a $\mathbb{Z}_4[x]$-submodule of $R$.
Hence, in this paper, we identify GQC codes of block lengths $(m_1,m_2,\ldots,m_l)$ over $\mathbb{Z}_4$ with $\mathbb{Z}_4[x]$-submodules of $R$.
\par
Define the canonical projection $\pi_i$ from $R$ to $R_i$ such that $$\pi_i((c_1(x)|\cdots|c_l(x)))=c_i(x).$$
If $\mathcal{C}$ is a GQC code of block lengths $(m_1,m_2,\ldots,m_l)$ over $\mathbb{Z}_4$, then $\pi_i(\mathcal{C})$ is a cyclic code over $\mathbb{Z}_4$ of length $m_i$ for $i=1,2,\ldots,l$.
\begin{proposition}
Let $\mathcal{C}$ be a GQC code of block lengths $(m_1,m_2,\ldots,m_l)$ over $\mathbb{Z}_4$. Then there exists a set $\{\mathbf{a}_1,\mathbf{a}_2,\ldots,\mathbf{a}_l\}$
such that $\mathcal{C}=\langle\mathbf{a}_1,\mathbf{a}_2,\ldots,\mathbf{a}_l\rangle$, where
\begin{align*}
&\mathbf{a}_1=(F_{1,1}|0|\cdots|0),\\
&\mathbf{a}_2=(F_{2,1}|F_{2,2}|\cdots|0),\\
&\vdots\\
&\mathbf{a}_l=(F_{l,1}|F_{l,2}|\cdots|F_{l,l}),
\end{align*}
 $F_{i,j}\in \mathbb{Z}_4[x]$ and $F_{i,i}=f_{i,i}(x)+2g_{i,i}(x)$, $f_{i,i}(x)$, $g_{i,i}(x)$ are monic polynomials over $\mathbb{Z}_4$ with $g_{i,i}(x)|f_{i,i}(x)|x^{m_i}-1$ for any $j=1,2,\ldots,i$, $i=1,2,\ldots,l$. The set $\{\mathbf{a}_1,\mathbf{a}_2,\ldots,\mathbf{a}_l\}$ is called the normalized generating set of $\mathcal{C}$.
\end{proposition}

\begin{proposition}
Let $\mathcal{C}$ be a GQC code of block lengths $(m_1,m_2,\ldots,m_l)$ over $\mathbb{Z}_4$ generated by the set $\{\mathbf{a}_1,\mathbf{a}_2,\ldots,\mathbf{a}_l\}$ given in Proposition 1. If nonzero elements in $\mathbf{a}_i$ are monic, then, for any $i=1,2,\ldots,l-1$, we can assume that ${\rm deg}(F_{i,i})> {\rm deg}(F_{i+1,i})>\cdots>{\rm deg}(F_{l,i})$.
\end{proposition}
\section{Minimum generating set of GQC codes}
 Let $\mathcal{C}$ be a linear code generated by the set $S$. If any element in $S$ can not be expressed as a linear combination of the others, then $S$ is called the minimum generating set of $\mathcal{C}$.
 In this section, we will give the minimum generating set of GQC codes of block lengths $(m_1,m_2,\ldots,m_l)$ over $\mathbb{Z}_4$.
\begin{proposition}
Let $\mathcal{C}$ be a GQC code of block lengths $(m_1,m_2,\ldots,m_l)$ over $\mathbb{Z}_4$ generated by the set $\{\mathbf{a}_1,\mathbf{a}_2,\ldots, \mathbf{a}_l\}$ given in Proposition 1.
Let ${\rm deg}(f_{i,i}(x))=t_i$, ${\rm deg}(g_{i,i}(x))=k_i$, $h_i(x)=\frac{x^{m_i}-1}{f_{i,i}(x)}$ for any $i=1,2,\ldots,l$. For $i=1,2,\ldots,l$, define
\begin{equation*}
S_{i,1}=\bigcup_{u=0}^{m_i-t_i-1}\left\{x^u\ast (F_{i,1}|F_{i,2}|\cdots|F_{i,i}|\cdots|0)\right\}
\end{equation*}
and
\begin{equation*}
S_{i,2}=\bigcup_{u=0}^{t_i-k_i-1}\left\{x^u\ast (h_i(x)F_{i,1}|h_i(x)F_{i,2}|\cdots|2h_i(x)g_{i,i}(x)|\cdots|0)\right\}.
\end{equation*}
Then $\bigcup_{i=1}^l{\left(S_{i,1}\bigcup S_{i,2}\right)}$ forms the minimum generating set of $\mathcal{C}$ as a $\mathbb{Z}_4[x]$-submodule of $R$.
\end{proposition}
\section{Dual codes of GQC codes}
In this section, we will investigate the dual codes of GQC codes and determine the relationship between GQC codes and their dual codes.
\par
Let $\mathcal{C}=\langle \mathbf{a}_1,\mathbf{a}_2,\ldots,\mathbf{a}_l \rangle$ be a GQC code of block lengths $(m_1,m_2,\ldots,m_l)$ over $\mathbb{Z}_4$.
Define the dual code $\mathcal{C}^\perp$ of $\mathcal{C}$ as
\begin{align*}
\mathcal{C}^\perp=\{&\mathbf{d}=(d_{1,0},d_{1,1},\ldots,d_{1,m_1-1}|\cdots|d_{l,0},d_{l,1},\ldots,d_{l,m_l-1})\in \mathbb{Z}_4^{m_1}\times\cdots\times\mathbb{Z}_4^{m_l}\\
&|\mathbf{d}\cdot\mathbf{c}=0\,\ {\rm for}\,\ \mathbf{c}=(c_{1,0},c_{1,1},\ldots,c_{1,m_1-1}|\cdots|c_{l,0},c_{l,1},\ldots,c_{l,m_l-1})\in\mathcal{C}\},
\end{align*}
where the dot product is the usual Euclidean inner product.
Let $m={\rm lcm}(m_1,m_2,\ldots,m_l)$. Define
\begin{align*}
f^*(x)=x^{{\rm deg}(f(x))}f\left({\frac{1}{x}}\right),~
\widetilde{f}(x)=f^*(x)x^{m-{\rm deg}(f(x))-1},
\end{align*}
where $f(x)\in \mathbb{Z}_4[x]$.
For any $$\mathbf{c}(x)=(c_1(x)|c_2(x)|\cdots|c_l(x)),~\mathbf{d}(x)=(d_1(x)|d_2(x)|\cdots|d_l(x))\in R,$$ define
$$\mathbf{d}(x)\circ \mathbf{c}(x)=d_1(x)\widetilde{c}_1(x)\theta_1(x^{m_1})+d_2(x)\widetilde{c}_2(x)\theta_2(x^{m_2})+\cdots+d_l(x)\widetilde{c}_l(x)\theta_l(x^{m_l}),$$
where $\theta_i(x^{m_i})=\frac{x^m-1}{x^{m_i}-1}$.
\par
From the above discussion, we can get the following result directly.
\begin{lemma}
Let $\mathbf{c}$, $\mathbf{d}\in \mathbb{Z}_4^{m_1}\times\cdots\times\mathbb{Z}_4^{m_l}$ be any two elements, $\mathbf{c}(x)$ and $\mathbf{d}(x)\in R$ be two elements corresponding to $\mathbf{c}$ and $\mathbf{d}$, respectively. Then $\mathbf{d}$ is orthogonal to $\mathbf{c}$ and its cyclic shifts if and only if  $\mathbf{d}(x)\circ \mathbf{c}(x)=0$.
\end{lemma}
From Lemma 4, we can redefine the dual code of $\mathcal{C}$ as follows
$$\mathcal{C}^\perp=\{\mathbf{d}(x)\in R|~\mathbf{d}(x)\circ\mathbf{c}(x)=0\,\ {\rm for\,\ any}\,\ \mathbf{c}(x)\in \mathcal{C}\}.$$
\par
For any $\mathbf{a}=(f_1(x)|f_2(x)|\cdots|f_l(x))\in R$, let $\mathbf{\overline{a}}=(\overline{f}_1(x)|\overline{f}_2(x)|\cdots|\overline{f}_l(x))$ and
$_{(i)}{\mathbf{a}}=(f_1(x)|f_2(x)|\cdots|f_i(x))\in R_1\times\cdots\times R_i$
be the truncated vector of $\mathbf{a}$ for any $i=1,2,\ldots,l$.
Let $\mathcal{C}=\langle \mathbf{a}_1,\mathbf{a}_2,\ldots,\mathbf{a}_l\rangle$ be a GQC code of block lengths $(m_1,m_2,\ldots,m_l)$ over $\mathbb{Z}_4$.
Let $$_{(i)}\mathcal{C}=\langle_{(i)}\mathbf{a}_1,_{(i)}\mathbf{a}_2,\ldots,_{(i)}\mathbf{a}_l\rangle, \,\
_{(i)}\mathbf{\varepsilon}=\langle_{(i)}\mathbf{e}_1, _{(i)}\mathbf{e}_2,\ldots,_{(i)}\mathbf{e}_i\rangle.$$
Clearly, $_{(i)}\mathcal{C}$ is a GQC code of block lengths $(m_1,m_2,\ldots,m_i)$ over $\mathbb{Z}_4$ and $(_{(i)}\mathcal{C})^\perp=_{(i)}\mathbf{\varepsilon}$.
Let
$$_{(i)}\overline{\mathcal{C}}=\langle_{(i)}\overline{\mathbf{a}}_1,_{(i)}
\overline{\mathbf{a}}_2,\ldots,_{(i)}\overline{\mathbf{a}}_l\rangle,~
_{(i)}\overline{\mathbf{\varepsilon}}=\langle_{(i)}\overline{\mathbf{e}}_1, _{(i)}\overline{\mathbf{e}}_2,\ldots,_{(i)}\overline{\mathbf{e}}_i\rangle.$$
Clearly, $_{(i)}\overline{\mathcal{C}}$ is a GQC code of block lengths $(m_1,m_2,\ldots,m_i)$ over $\mathbb{Z}_2$ and $({_{(i)}\overline{\mathcal{C}}})^\perp=_{(i)}\overline{\mathcal{\varepsilon}}$ since the map $-$ is a surjective ring homomorphism.
\par
Let $\langle \mathbf{g}_{1},\mathbf{g}_2,\ldots,\mathbf{g}_i\rangle$ be the normalized generating set of ${_{(i)}\mathcal{C}}$, where
\begin{align*}
&\mathbf{g}_{1}=(G_{1,1}|0|\cdots|0),\\
&\mathbf{g}_{2}=(G_{2,1}|G_{2,2}|\cdots|0),\\
&\vdots\\
&\mathbf{g}_{i}=(G_{i,1}|G_{i,2}|\cdots|G_{i,i}),\\
\end{align*}
$G_{r,r}$ are monic with $\overline{G}_{r,r}$~no~multiple~roots for $r=1,2\ldots,i$ and $G_{r,r-1}$ are monic with $\overline{G}_{r,r-1}$~no~multiple~roots for $r=2\ldots,i$.
Then the normalized generating set of ${_{(i)}\overline{\mathcal{C}}}$ can be given as
\begin{align*}
&\overline{\mathbf{g}}_{1}=(\overline{G}_{1,1}|0|\cdots|0),\\
&\overline{\mathbf{g}}_{2}=(\overline{G}_{2,1}|\overline{G}_{2,2}|\cdots|0),\\
&\vdots\\
&\overline{\mathbf{g}}_{i}=(\overline{G}_{i,1}|\overline{G}_{i,2}|\cdots|\overline{G}_{i,i}).\\
\end{align*}
\par
For any $r=1,2,\ldots,i$, define
$A_{r,r}^{(0)}=G_{r,r}^{(0)}=G_{r,r},\,\  B_{r,r}^{(0)}=0, \,\  \mathbf{g}_r^{(0)}=\mathbf{g}_r.$
For any $r=2,\ldots,i$, let
\begin{align*}
&A_{r,r-1}^{(1)}=\frac{G_{r-1,r-1}^{(0)}}{{\rm gcd_4}\left(G_{r-1,r-1}^{(0)}, G_{r,r-1}^{(0)}\right)},\,\
B_{r,r-1}^{(1)}=\frac{G_{r,r-1}^{(0)}}{{\rm gcd_4}\left(G_{r-1,r-1}^{(0)}, G_{r,r-1}^{(0)}\right)},\\
&\mathbf{g}_r^{(1)}=A_{r,r-1}^{(1)}\mathbf{g}_r^{(0)}-B_{r,r-1}^{(1)}\mathbf{g}_{r-1}^{(0)}=(G_{r,1}^{(1)}|G_{r,2}^{(1)}|\cdots|G_{r,r}^{(1)}|0|\cdots|0).
\end{align*}
Clearly, $G_{r,r-1}^{(1)}=0$ and $G_{r,r}^{(1)}=A_{r,r}^{(0)}A_{r,r-1}^{(1)}$. For any $r=3,\ldots,i$, let $G_{r,r-2}^{(1)}$ be monic with $\overline{G}_{r,r-2}^{(1)}$~no~multiple~roots, and
\begin{align*}
&A_{r,r-2}^{(2)}=\frac{G_{r-2,r-2}^{(0)}}{{\rm gcd_4}\left(G_{r-2,r-2}^{(0)}, G_{r,r-2}^{(1)}\right)},\,\
B_{r,r-2}^{(2)}=\frac{G_{r,r-2}^{(1)}}{{\rm gcd_4}\left(G_{r-2,r-2}^{(0)}, G_{r,r-2}^{(1)}\right)},\\
&\mathbf{g}_r^{(2)}=A_{r,r-2}^{(2)}\mathbf{g}_r^{(1)}-B_{r,r-2}^{(2)}\mathbf{g}_{r-2}^{(0)}=(G_{r,1}^{(2)}|G_{r,2}^{(2)}|\cdots|G_{r,r}^{(2)}|0|\cdots|0).
\end{align*}
Clearly, $G_{r,r-1}^{(2)}=G_{r,r-2}^{(2)}=0$ and $G_{r,r}^{(2)}=A_{r,r}^{(0)}A_{r,r-1}^{(1)}A_{r,r-2}^{(2)}$. For any $r=4,\ldots,i$, if $G_{r,r-3}^{(2)}$ are monic with $\overline{G}_{r,r-3}^{(2)}$~no~multiple~roots, then repeat the above process. Hence, for any $k=0,1,\ldots,r-2$, $r=2,3,\ldots,i$, if $G_{r,r-k-1}^{(k)}$ are monic with $\overline{G}_{r,r-k-1}^{(k)}$~no~multiple~roots, then we have
\begin{align*}
\mathbf{g}_{r}^{(r-1)}=(0|0|\cdots|G_{r,r}^{(r-1)}|0\cdots|0),
\end{align*}
where $G_{r,r}^{(r-1)}=A_{r,r}^{(0)}A_{r,r-1}^{(1)}\cdots A_{r,1}^{(r-1)}$.
\par
For the normalized generating set of $_{(i)}\overline{\mathcal{C}}$, we can get it by the above method as follows
\begin{align*}
&\overline{A}_{r,r}^{(0)}=\overline{G}_{r,r}^{(0)}=\overline{G}_{r,r},\,\  \overline{B}_{r,r}^{(0)}=0, \,\  \overline{\mathbf{g}}_r^{(0)}=\overline{\mathbf{g}}_r,\\
&\overline{A}_{r,r-k}^{(k)}=\frac{\overline{G}_{r-k,r-k}^{(0)}}{{\rm gcd}\left(\overline{G}_{r-k,r-k}^{(0)}, \overline{G}_{r,r-k}^{(k-1)}\right)},\,\
\overline{B}_{r,r-k}^{(k)}=\frac{\overline{G}_{r,r-k}^{(k-1)}}{{\rm gcd}\left(\overline{G}_{r-k,r-k}^{(0)},\overline{G}_{r,r-k}^{(k-1)}\right)},\\
&\overline{\mathbf{g}}_r^{(k)}=\overline{A}_{r,r-k}^{(k)}\overline{\mathbf{g}}_r^{(k-1)}+\overline{B}_{r,r-k}^{(k)}\overline{\mathbf{g}}_{r-k}^{(0)}
=\left(\overline{G}_{r,1}^{(k)}|\overline{G}_{r,2}^{(k)}|\cdots|\overline{G}_{r,r}^{(k)}|0|\cdots|0\right),
\end{align*}
where $\overline{G}_{r,r}^{(r-1)}=\overline{A}_{r,r}^{(0)}\overline{A}_{r,r-1}^{(1)}\cdots \overline{A}_{r,1}^{(r-1)}$ and
$\overline{\mathbf{g}}_{r}^{(r-1)}=\left(0|0|\cdots|\overline{G}_{r,r}^{(r-1)}|0\cdots|0\right).$
\par
In the rest of this section, suppose that $G_{r,r}$ are monic with ${\rm gcd}(\overline{G}_{r,r},\overline{G}'_{r,r})=1$ for any $r=1,2,\ldots,i$, $i=1,2,\ldots,l$.
Let $\{\mathbf{e}_1,\mathbf{e}_2,\ldots,\mathbf{e}_l\}$ be the normalized generating set of $\mathcal{C}^\perp$, where
\begin{align*}
&\mathbf{e}_1=(E_{1,1}|0|\cdots|0),\\
&\mathbf{e}_2=(E_{2,1}|E_{2,2}|\cdots|0),\\
&\vdots\\
&\mathbf{e}_l=(E_{l,1}|E_{l,2}|\cdots|E_{l,l}),
\end{align*}
$E_{i,j}\in\mathbb{Z}_4[x]$ and $E_{i,i}=u_{i,i}(x)+2v_{i,i}(x)$, $u_{i,i}(x)$ and $v_{i,i}(x)$ are monic polynomials over $\mathbb{Z}_4$ with $v_{i,i}(x)|u_{i,i}(x)|x^{m_i}-1$ for $j=1,2,\ldots,i$, $i=1,2,\ldots,l$.

\begin{lemma}
Let $\mathcal{C}$ be a GQC code of block lengths $(m_1,m_2,\ldots,m_l)$ over $\mathbb{Z}_4$ generated by the set $\{\mathbf{a}_1,\mathbf{a}_2,\ldots,\mathbf{a}_l\}$ given in Proposition 1 and
$\mathcal{C}^\perp=\langle\mathbf{e}_1,\mathbf{e}_2,\ldots,\mathbf{e}_l\rangle$ be its dual code.
Let
$_{(i)}\mathcal{C}=\langle \mathbf{g}_1,\mathbf{g}_2,\ldots,\mathbf{g}_i\rangle$ be the truncated code associated to $\mathcal{C}$. For any $k=0,1,\ldots,r-2$, $r=2,\ldots,i$ if $G_{r,r-k-1}^{(k)}$ are monic with $\overline{G}_{r,r-k-1}^{(k)}$~no~multiple~roots, then
\begin{align*}
E_{t,r}{G_{r,r}^{(r-1)}}^*\equiv 0 ({\rm mod} \,\ x^{m_r}-1),
\end{align*}
where $r=1,2,\ldots,i$ and $t=r,r+1,\ldots,i$.
\end{lemma}
\par
Particularly, if $t=i$, then $E_{i,r}{G_{r,r}^{(r-1)}}^*\equiv 0~({\rm mod} ~ x^{m_r}-1)$ for any $r=1,2\ldots,i$, which implies that there exists $\lambda_{i,r}\in \mathbb{Z}_4[x]$
such that $E_{i,r}{G_{r,r}^{(r-1)}}^*=(x^{m_r}-1)\lambda_{i,r}$. In the following, we will determine the degree of $E_{i,i}$ by the Hensel lift given in \cite{Wan}.
\begin{lemma}\cite[Proposition 5.12]{Wan}(Hensel lift)
Let $n$ be an odd positive integer and $f_1(x)|x^{n}-1$ be a polynomial over $\mathbb{Z}_2$. Then there exists a unique monic polynomial $f(x)|x^{n}-1$ over $\mathbb{Z}_4$ and $\overline{f}(x)=f_1(x)$.
\end{lemma}
\begin{lemma}
Let $\mathcal{C}$ be a GQC code of block lengths $(m_1,m_2,\ldots,m_l)$ over $\mathbb{Z}_4$ generated by the set $\{\mathbf{a}_1,\mathbf{a}_2,\ldots,\mathbf{a}_l\}$ given in Proposition 1 and
$\mathcal{C}^\perp=\langle\mathbf{e}_1,\mathbf{e}_2,\ldots,\mathbf{e}_l\rangle$ be its dual code. Let
$_{(i)}\mathcal{C}=\langle \mathbf{g}_1,\mathbf{g}_2,\ldots,\mathbf{g}_i\rangle$ be the truncated code associated to $\mathcal{C}$.
If $F_{i,i}|x^{m_i}-1$ for $i=1,2,\ldots,l$ and  $G_{i,i-k-1}^{(k)}$ are monic with $\overline{G}_{i,i-k-1}^{(k)}$~no~multiple~roots for $k=0,1,\ldots,i-2$, $i=2,\ldots,l$, then
\begin{align*}
{\rm deg}(E_{i,i})=m_i-{\rm deg}(A_{i,i}^{(0)})-{\rm deg}(A_{i,i-1}^{(1)})-\cdots-{\rm deg}(A_{i,1}^{(i-1)}),
\end{align*}
for any $i=1,2,\ldots,l$.
\end{lemma}
\par
By Lemmas 5 and 7, we can get the following proposition.
\begin{proposition}
Let $\mathcal{C}$ be a GQC code of block lengths $(m_1,m_2,\ldots,m_l)$ over $\mathbb{Z}_4$ generated by the set $\{\mathbf{a}_1,\mathbf{a}_2,\ldots,\mathbf{a}_l\}$ given in Proposition 1 and
$\mathcal{C}^\perp=\langle\mathbf{e}_1,\mathbf{e}_2,\ldots,\mathbf{e}_l\rangle$ be its dual code.
Let
$_{(i)}\mathcal{C}=\langle \mathbf{g}_1,\mathbf{g}_2,\ldots,\mathbf{g}_i\rangle$ be the truncated code associated to $\mathcal{C}$. If $F_{i,i}|x^{m_i}-1$ for $i=1,2,\ldots,l$ and  $G_{i,i-k-1}^{(k)}$ are monic with $\overline{G}_{i,i-k-1}^{(k)}$~no~multiple~roots for $k=0,1,\ldots,i-2$, $i=2,\ldots,l$, then
\begin{align*}
E_{i,i}{G_{i,i}^{(i-1)}}^*=x^{m_i}-1~{\rm or}~1-x^{m_i},
\end{align*}
for any $i=1,2,\ldots,l$.
\end{proposition}
\par
In the following, we will determine $E_{i,r}$, where $r=1,2,\ldots,i-1$.
\par
\begin{proposition}
Let $\mathcal{C}$ be a GQC code of block lengths $(m_1,m_2,\ldots,m_l)$ over $\mathbb{Z}_4$ generated by the set $\{\mathbf{a}_1,\mathbf{a}_2,\ldots,\mathbf{a}_l\}$ given in Proposition 1, and
$\mathcal{C}^\perp=\langle\mathbf{e}_1,\mathbf{e}_2,\ldots,\mathbf{e}_l\rangle$ be its dual code. Let
$_{(i)}\mathcal{C}=\langle \mathbf{g}_1,\mathbf{g}_2,\ldots,\mathbf{g}_i\rangle$ be the truncated code associated to $\mathcal{C}$.
If $F_{i,i}|x^{m_i}-1$ for $i=1,2,\ldots,l$ and  $G_{r,r-k-1}^{(k)}$ are monic with $\overline{G}_{r,r-k-1}^{(k)}$~no~multiple~roots for $k=0,1,\ldots,r-2$, $r=2,\ldots,i-1$, then
$$E_{i,r}{G_{r,r}^{(r-1)}}^*=(x^{m_r}-1)\lambda_{i,r},$$
where
\begin{align*}
\lambda_{i,1}&\equiv3\lambda_{i,i}\left({B_{i,1}^{(i-1)}}^*\right)^{-1}x^{m'+{\rm deg}(G_{i,1}^{(i-2)})-{\rm deg}(G_{i,i}^{(i-2)})}~({{\rm mod}~A_{i,1}^{(i-1)}}^*),\\
\lambda_{i,s}&\equiv 3\lambda_{i,1}{B_{s,1}^{(s-1)}}^*x^{m'+{\rm deg}(G_{s,s}^{(s-2)})-{\rm deg}(G_{s,1}^{(s-2)})}~({\rm mod}~{A_{s,1}^{(s-1)}}^*),
\end{align*}
and $m'={\rm lcm}(m_1,m_2,\ldots,m_i)$ for $s=2\ldots,i-1$, $i=1,2,\ldots,l$.
\end{proposition}
\section{Conclusion}
In this paper, we studied the normalized generating set and the minimum generating set of GQC codes over $\mathbb{Z}_4$. The results in this paper are also valid on GQC codes over finite chain rings with nilpotency index 2. It is an interesting open problem to study GQC codes over finite chain rings with nilpotency index $\geq 3$. Another interesting open problem is to consider quantum codes construction from GQC codes over $\mathbb{Z}_4$.

\bibliographystyle{model1a-num-names}

\end{document}